\documentclass[fleqn]{mat01}
\usepackage{times,mathtimy,amssymb,latexsym}
\begin{document}

\def\d{\mbox{\rm d}}

\newtheorem{theore}{Theorem}
\renewcommand\thetheore{\arabic{section}.\arabic{theore}}
\newtheorem{theor}[theore]{\bf Theorem}
\newtheorem{lem}[theore]{Lemma}
\newtheorem{rem}[theore]{Remark}
\newtheorem{coro}[theore]{\rm COROLLARY}

\newtheorem{step}{Step}

\renewcommand\theequation{\arabic{section}.\arabic{equation}}

\setcounter{page}{241}
\firstpage{241}

\title{Some properties of complex matrix-variate generalized Dirichlet
integrals}

\markboth{Joy Jacob, Sebastian George and A~M~Mathai}{Complex
matrix-variate beta random variables}

\author{JOY JACOB, SEBASTIAN GEORGE and A~M~MATHAI$^{*}$}

\address{Department of Statistics, St.~Thomas College, Arunapuram P.O.,
Palai, Kottayam~686~574, India\\
\noindent $^{*}$McGill University, Montreal, Canada and
Centre for Mathematical Sciences, Pala Campus, Arunapuram P.O., Pala 686 574, India\\
\noindent E-mail: jjstc@sancharnet.in}

\volume{115}

\mon{August}

\parts{3}

\pubyear{2005}

\Date{MS received 9 October 2003; revised 14 October 2004}

\begin{abstract}
Dirichlet integrals and the associated Dirichlet statistical densities
are widely used in various areas. Generalizations of Dirichlet integrals
and Dirichlet models to matrix-variate cases, when the matrices are real
symmetric positive definite or hermitian positive definite, are
available \cite{4}. Real scalar variables case of the
Dirichlet models are generalized in various directions. One such
generalization of the type-2 or inverted Dirichlet is looked into in
this article. Matrix-variate analogue, when the matrices are hermitian
positive definite, are worked out along with some properties which are
mathematically and statistically interesting.
\end{abstract}


\maketitle

\section{Introduction}

This paper deals with probability densities on the space of matrices.
All the matrices appearing in this article are $p \times p$ hermitian
positive definite unless specified otherwise. $X, Y, \ldots$ will denote
matrices whose elements are functionally independent real scalar
mathematical variables or random variables. $\tilde{X}, \tilde{Y}, \ldots$
will denote matrices whose elements are in the complex domain. Constant
matrices will be denoted by $A, B, \ldots$ whether the elements are real or
complex. ${\rm tr}(\cdot)$ will denote the trace of the matrix
$(\cdot), |(\cdot)|$ will denote the determinant as well as absolute
value of $(\cdot)$ and $|\det (\cdot)|$ will denote the absolute
value of the determinant of $(\cdot)$. Transpose will be denoted by a
prime, complex conjugate by a bar, conjugate transpose by a star. Thus
$\tilde{X} = \tilde{X}^* > 0$ will mean the hermitian matrix $\tilde{X}$ is
positive definite. $\d X$ will indicate the wedge product of
differentials in $X$. For example, when $X = (x_{ij})$ and real,
\begin{align*}
\d X &= \wedge_{i,j} \d x_{ij} \quad \hbox{when all} \ x_{ij}\hbox{'s} \
\hbox{are distinct}\\[.2pc]
&= \wedge_{i \geq j} \d x_{ij} = \wedge_{i \leq j} \d x_{ij} \quad \hbox{when} \ X = X'.
\end{align*}
$\int_{0 < \tilde{X} = \tilde{X}^* < I}(\cdot) \d \tilde{X}$ indicates the
integral of $(\cdot)$ over all $\tilde{X}$ such that
$\tilde{X} = \tilde{X}^* > 0$ as well as $I - \tilde{X} > 0$. In other words,
all eigenvalues of $\tilde{X}$ are in the open interval (0,1), where $I$
denotes an identity matrix. For any positive definite hermitian matrix
$A$ we denote by $A^{\frac{1}{2}}$ the hermitian positive definite square
root of $A$.

In our discussions, we need a few Jacobians of matrix
transformations and integrals over real scalar functions of matrix
arguments. These will be stated here without proof. For proofs and other
details, see \cite{4}.
\begin{equation}
\tilde{Y} = A \tilde{X} A^{*}, |A| \neq 0 \Rightarrow \d \tilde{Y}
= |\det (A)|^{2p} \d \tilde{X} = |\det (AA^{*})|^p \d \tilde{X},
\end{equation}
where $A$ is a constant matrix (p.~183 of \cite{4})
\begin{equation}
\tilde{Y} = \tilde{X}^{-1},~|\tilde{X}| \neq 0 \Rightarrow \d
\tilde{Y} = |\det (\tilde{X})|^{-2p} \d \tilde{X} \quad \hbox{for} \ \tilde{X} =
\tilde{X}^*,
\end{equation}
(p.~190 of \cite{4}). Let $\tilde{T} =
(\tilde{t}_{ij}),~\tilde{t}_{ij} = 0$ for $i < j,~\tilde{t}_{jj} =
t_{jj} > 0$ (real and positive) for $j = 1, \ldots, p$. That is,
$\tilde{T}$ is a lower triangular matrix with real positive
diagonal elements. Then
\begin{equation}
\tilde{X} = \tilde{T}\tilde{T}^{*} \Rightarrow \d \tilde{X} = 2^p
\left\{ \prod_{j=1}^p t_{jj}^{2(p-j) + 1}\right\} \d \tilde{T}.
\end{equation}
The complex matrix-variate gamma, denoted by
$\tilde{\Gamma}_p(\alpha)$, and the gamma integral are defined as
follows:
\begin{align}
\tilde{\Gamma}_p(\alpha) &= \pi^{{\frac{p(p-1)}{2}}} \Gamma(\alpha)
\Gamma(\alpha-1) \cdots \Gamma(\alpha - p + 1),~\Re (\alpha) > p - 1\\[.2pc]
&= \int_{\tilde{X} = {\tilde{X}}^{*} > 0} |\det
(\tilde{X})|^{\alpha - p}{\rm e}^ {-{\rm tr}(\tilde{X})} \d
\tilde{X},
\end{align}
where $\Re(\cdot)$ denotes
the real part of $(\cdot)$, (p.~188 of \cite{4}). The integral in (1.5)
when evaluated with the help of (1.3) yields (1.4). The complex
matrix-variate beta, denoted by $\tilde{B}_p(\alpha, \beta)$, and the
beta integrals are defined as follows (p.~198 of \cite{4}):
\begin{align}
\tilde{B}_p(\alpha, \beta) &= {\frac{\tilde{\Gamma}_p(\alpha)
\tilde{\Gamma}_p(\beta)}{\tilde{\Gamma}_p(\alpha+\beta)}},
~\Re (\alpha) > p - 1, \quad \Re (\beta) > p - 1\\[.2pc]
&= \int_{0 < \tilde{X} = \tilde{X}^* < I}|\det (\tilde{X})|^{\alpha - p}
|\det (I - \tilde{X})|^{\beta - p} \d \tilde{X}\\[.2pc]
&= \int_{\tilde{X} = \tilde{X}^* > 0}|\det (\tilde{X})|^{\alpha - p}
|\det (I + \tilde{X})|^{-(\alpha + \beta)} \d \tilde{X}.
\end{align}
The functions of $\tilde{X}$ appearing in (1.7) and (1.8), divided by
the normalizing constant $\tilde{B}_p (\alpha, \beta)$, produce the
complex matrix-variate type-1 and type-2 beta densities respectively
(p.~357 of \cite{4}). Dirichlet integrals and the associated Dirichlet
densities come in naturally in order statistics problems, in reliability
analysis and in certain survival analysis problems. In Bayesian
statistical analysis a Dirichlet model is usually taken as the prior
distribution for multinomial probabilities. In random division and
certain geometrical probability problems Dirichlet model comes in
naturally, see for example \cite{5}. Dirichlet model is extended to
matrix-variate Liouville type by \cite{2}. See \cite{1} for some
applications in engineering, \cite{3} for statistical
applications and \cite{6} for some physics problems.

The model that we are going to deal with in the present article is the
complex matrix-variate analogue of an extended inverted or type-2
Dirichlet model incorporating successive sums of the variables into it.

\section{Matrix-variate analogue of an extended Dirichlet model}

Let $\tilde{X_1}, \ldots, \tilde{X_k}$ be hermitian positive definite
$p\times p$ matrix random variables having the joint density function
\setcounter{equation}{0}
\begin{align}
f(\tilde{X}_1, \ldots, \tilde{X}_k) &= c_k|\det (\tilde{X}_1)|^
{\alpha_1 - p} \cdots |\det (\tilde{X}_k)|^{\alpha_k - p}\nonumber\\[.2pc]
&\quad \, \times|\det (I + \tilde{X}_2 + \cdots + \tilde{X}_k)|^{\beta_1}\cdots
|\det (I + \tilde{X}_k)|^{\beta_{k - 1}}\nonumber\\[.2pc]
&\quad \, \times|\det (I + \tilde{X}_1 + \cdots + \tilde{X}_k)|^
{-(\alpha_1 + \cdots + \alpha_{k + 1} + \beta_1 + \cdots + \beta_k)}
\end{align}
for
$\Re (\alpha_j) > p - 1,~\Re(\alpha_{j+1} + \cdots + \alpha_{k + 1}
+ \beta_j + \cdots + \beta_k) > p-1,~j = 1, \ldots, k$, and
$f(\tilde{X_1}, \ldots, \tilde{X_k}) = 0$ elsewhere, where $c_k$ is the
normalizing constant. We will study some properties of (2.1) in the
present article. Since
\begin{equation*}
\int_{\tilde{X_1}} \cdots \int_{\tilde{X_k}} f (\tilde{X_1},
\ldots, \tilde{X_k}) \d \tilde{X_1} \wedge \ldots \wedge \d
\tilde{X_k} = 1,
\end{equation*}
we can evaluate $c_k$ by successive integration starting with
$\tilde{X_1}$. For fixed $\tilde{X}_2, \ldots, \tilde{X}_k$ let
\begin{align*}
L_1 &= \int_{\tilde{X}_1 = \tilde{X}_1^{*} > 0}|\det (\tilde{X}_1)|^
{\alpha_1-p}\\[.2pc]
&\quad \times |\det (I + \tilde{X}_1 + \cdots + \tilde{X}_k)|^
{-(\alpha_1 + \cdots + \alpha_{k + 1} + \beta_1 + \cdots + \beta_k)} \d
\tilde{X}_1.
\end{align*}
Note that
\begin{align*}
|\det (I + \tilde{X}_1 + \cdots + \tilde{X}_k)| &= |\det (I +
\tilde{X}_2 + \cdots + \tilde{X}_k)|\\[.2pc]
&\quad \, \times|\det [I + (I + \tilde{X}_2 + \cdots +
\tilde{X}_k)^{-{1/2}}\\[.2pc]
&\quad \, \times \tilde{X}_1(I + \tilde{X}_2 + \cdots +
\tilde{X}_k)^{-{{1}/{2}}}]|.
\end{align*}
Now, make the transformation
\begin{equation*}
\tilde{Y}_1 = (I + \tilde{X}_2 + \cdots +
\tilde{X}_k)^{-{{1}/{2}}} \tilde{X}_1 (I + \tilde{X}_2 + \cdots +
\tilde{X}_k)^{-{{1}/{2}}}.
\end{equation*}
Then
\begin{equation*}
\d \tilde{Y}_1 = |\det (I + \tilde{X}_2 + \cdots + \tilde{X}_k)^{-p}| \d
\tilde{X}_1.
\end{equation*}
Substituting for $\tilde{X}_1$ in terms of $\tilde{Y}_1$ we have
\begin{align*}
L_1 &= |\det (I + \tilde{X}_2 + \cdots + \tilde{X}_k)|^{-(\alpha_2 +
 \cdots + \alpha_{k + 1} + \beta_1 + \cdots + \beta_k)}\\[.2pc]
&\quad \, \times\int_{\tilde{Y}_1 = \tilde{Y}_1^{*} > 0}|\det (\tilde{Y}_1)|^
{\alpha_1 - p}|\det (I + \tilde{Y}_1)|^ {-(\alpha_1 + \cdots + \alpha_{k + 1} +
\beta_1 + \cdots + \beta_k)} \d \tilde{Y}_1\\[.2pc]
&=|\det (I + \tilde{X}_2 + \cdots + \tilde{X}_k)|^{-(\alpha_2 + \cdots +
\alpha_{k + 1} + \beta_1 + \cdots + \beta_k)}\\[.2pc]
&\quad \, \times {\frac{\tilde{\Gamma}_p(\alpha_1) \tilde{\Gamma}_p (\alpha_2 + \cdots +
\alpha_{k + 1} + \beta_1 + \cdots + \beta_k)}{\tilde{\Gamma}_p(\alpha_1 +
\cdots + \alpha_{k + 1} + \beta_1 + \cdots + \beta_k)}}
\end{align*}
for $\Re(\alpha_1) > p - 1, ~\Re(\alpha_2 + \cdots + \alpha_{k + 1} + \beta_1 +
\cdots + \beta_k) > p - 1$. The $\tilde{Y}_1$-integral is evaluated by using
the type-2 beta integral of (1.8). Successive integrations of
$\tilde{X}_2, \ldots, \tilde{X}_k$ will yield the result as follows:
\begin{equation}
c_k^{-1} = \prod_{j = 1}^k {\frac{\tilde{\Gamma}_p(\alpha_j)\tilde{\Gamma}_p(\alpha_{j + 1} +\cdots +
 \alpha_{k + 1} + \beta_j + \cdots + \beta_k)}{\tilde{\Gamma}_p(\alpha_j +
 \cdots + \alpha_{k + 1} + \beta_j + \cdots + \beta_k)}}
\end{equation}
for $\Re(\alpha_j) > p - 1,~\Re(\alpha_{j + 1} + \cdots + \alpha_{k + 1}
+ \beta_j + \cdots + \beta_k) > p - 1, ~j = 1, \ldots, k$.

We can look into some interesting results and the corresponding matrix
transformations. These will be stated as theorems. We need two more
results in order to establish our main results. These will be listed as
lemmas.

\begin{lem}
Let $\tilde{X}$ and $A$ be $p \times p$ hermitian positive definite
matrices. Then $\tilde{Y}_1 = A^{{1}/{2}} \tilde{X} A^{{1}/{2}}$ and
$\tilde{Y}_2 = \tilde{X}^{{1}/{2}} A\tilde{X}^{{1}/{2}}$ have the same
eigenvalues.

This can be easily seen by looking at the determinantal equations for
the eigenvalues
\begin{equation*}
|A^{{1}/{2}}\tilde{X}A^{{1}/{2}} - \lambda I| = 0 \Rightarrow
 |\tilde{X} A - \lambda I| = 0 \Rightarrow| \tilde{X}^{{1}/{2}} A \tilde{X}^
{{1}/{2}} - \lambda I| = 0.
\end{equation*}
Thus{\rm ,} $\tilde{Y}_1$ and $\tilde{Y}_2$ have the same equations giving
rise to the same eigenvalues $\lambda_1, \ldots,
\lambda_p$, $\lambda_j > 0,~~j = 1, \ldots, p$.
\end{lem}

\begin{lem}
Let the common real eigenvalues of $\tilde{Y}_1$ and $\tilde{Y}_2$ of
Lemma~$2.1$ be distinct. Then the wedge product $\d \tilde{Y}_1 = \d
\tilde{Y}_2$ in the integrals.
\end{lem}

\begin{proof}
Let $\tilde{U}$ and $\tilde{V}$ be unitary matrices with real diagonal
elements such that
\begin{equation*}
\tilde{W}_1 = \tilde{U}^*\tilde{Y}_1\tilde{U} = D = {\rm diag}(\lambda_1, \ldots, \lambda_p)=
 \tilde{V}^*\tilde{Y}_2\tilde{V} = \tilde{W}_2.
\end{equation*}
Then from Theorem~4.4 of \cite{4}
\begin{equation*}
\d \tilde{Y}_1 = \d \tilde{W}_1 = \left\{\prod_{j > k}
  |\lambda_k - \lambda_j|^2\right\} \d D\wedge \d \tilde{G}_1
\end{equation*}
and
\begin{equation*}
\d \tilde{Y}_2 = \d \tilde{W}_2 = \left\{\prod_{j > k}
|\lambda_k - \lambda_j|^2\right\} \d D \wedge \d \tilde{G}_2,
\end{equation*}
where $\d \tilde{G}_1$ and $\d\tilde{G}_2$ are the following:
\begin{equation*}
\d \tilde{G}_1 = \wedge[\tilde{U} (\d \tilde{U})] \quad \hbox{and} \quad
\d \tilde{G}_2 = \wedge[\tilde{V} (\d \tilde{V})].
\end{equation*}
$(\d \tilde{U})$ and $(\d \tilde{V})$ denote the matrices of
differentials (entry-wise derivatives) in $\tilde{U}$ and $\tilde{V}$
respectively. Now from Corollary 4.3.1 of \cite{4}
\begin{equation*}
\int_{\tilde{O}_1(p)} \d \tilde{G}_1 = \int_{\tilde{O}_1(p)} \d
\tilde{G}_2 ={\frac{\pi^{p(p-1)}}{\tilde{\Gamma}_p(p)}},
\end{equation*}
where $\tilde{O}_1(p)$ is the orthogonal group of unitary matrices with
real diagonal elements.
\end{proof}

\setcounter{theore}{0}
\begin{theor}[\!]
Let $\tilde{X}_1, \ldots, \tilde{X}_k$ be matrix-variate random
variables having the joint distribution as in $(2.1)$. Consider
the transformation
\begin{align}
\tilde{Y}_1 &= (I + \tilde{X}_1 + \cdots + \tilde{X}_k)^{-{{1}/{2}}}
 \tilde{X}_1(I + \tilde{X}_1 + \cdots + \tilde{X}_k)^{-{{1}/{2}}}\nonumber\\[.2pc]
 \tilde{Y}_2 &= (I + \tilde{X}_2 + \cdots + \tilde{X}_k)^{-{{1}/{2}}}
 \tilde{X}_2(I + \tilde{X}_2 + \cdots + \tilde{X}_k)^{-{{1}/{2}}}\nonumber\\[.2pc]
 & \ \ \vdots\nonumber\\[.2pc]
 \tilde{Y}_k &= (I + \tilde{X}_k)^{-{{1}/{2}}}\tilde{X}_k (I + \tilde{X}_k)^{-{{1}/{2}}}.
\end{align}
Then $\tilde{Y}_1, \ldots, \tilde{Y}_k$ are independent{\rm ,} and
further{\rm ,} $\tilde{Y}_j$ has a type-$1$ beta density with the
parameters $(\alpha_j,~\alpha_{j+1} + \cdots + \alpha_{k + 1} +
\beta_j + \cdots + \beta_k),~j = 1, \ldots, k$.
\end{theor}

\begin{proof}
From the transformation in (2.3), we have
\begin{align}
I-\tilde{Y}_1 &= I-(I + \tilde{X}_1 + \cdots + \tilde{X}_k)^{-{{1}/{2}}}
  \tilde{X}_1(I + \tilde{X}_1 + \cdots + \tilde{X}_k)^{-{{1}/{2}}}\nonumber\\[.2pc]
  &= (I + \tilde{X}_1 + \cdots + \tilde{X}_k)^{-{{1}/{2}}}[(I + \tilde{X}_1
  + \cdots + \tilde{X}_k) - \tilde{X}_1]\nonumber\\[.2pc]
 &\quad \, \times (I + \tilde{X}_1 + \cdots + \tilde{X}_k)^{-
{{1}/{2}}}\nonumber\\[.2pc]
  &= (I + \tilde{X}_1 + \cdots + \tilde{X}_k)^{-{{1}/{2}}}(I + \tilde{X}_2 + \cdots
  + \tilde{X}_k)\nonumber\\[.2pc]
&\quad \, \times (I + \tilde{X}_1 + \cdots + \tilde{X}_k)^{-{{1}/{2}}}\nonumber\\[.2pc]
  I-\tilde{Y}_2 &= (I + \tilde{X}_2 + \cdots + \tilde{X}_k)^{-{{1}/{2}}} (I + \tilde{X}_3 + \cdots
  + \tilde{X}_k)\nonumber\\[.2pc]
&\quad \, \times (I + \tilde{X}_2 + \cdots + \tilde{X}_k)^{-{{1}/{2}}}\nonumber\\[.2pc]
  I-\tilde{Y}_{k-1} &= (I + \tilde{X}_{k - 1} + \tilde{X}_k)^{-{{1}/{2}}}(I + \tilde{X}_k)
 (I + \tilde{X}_{k - 1} + \tilde{X}_k)^{-{{1}/{2}}}\nonumber\\[.2pc]
  I-\tilde{Y}_k &= (I + \tilde{X}_k)^{-1}.
\end{align}
From the above representations of $I-\tilde{Y}_1, \ldots,
I-\tilde{Y}_k$, from (1.1), (1.2) and from Lemma~2.2, we can evaulate the
Jacobian of the transformation in (2.3) from (2.4).
\begin{align*}
\d \tilde{Y}_1 &= \d [(I + \tilde{X}_1 + \cdots + \tilde{X}_k
)^{-{{1}/{2}}}(I + \tilde{X}_2 + \cdots + \tilde{X}_k)\\[.2pc]
&\quad \, \times (I + \tilde{X}_1 + \cdots + \tilde{X}_k
)^{-{{1}/{2}}}]\\[.2pc]
&= \d [(I + \tilde{X}_2 + \cdots + \tilde{X}_k)^{{1}/{2}}(I + \tilde{X}_1 + \cdots
+ \tilde{X}_k)^{-1}\\[.2pc]
&\quad \, \times (I + \tilde{X}_2 + \cdots + \tilde{X}_k)^{{1}/{2}}]\\[.2pc]
&= |\det (I + \tilde{X}_2 + \cdots + \tilde{X}_k)|^p|\det (I + \tilde{X}_1 + \cdots
+ \tilde{X}_k)|^{-2p} \d \tilde{X}_1
\end{align*}
for fixed $\tilde{X}_2, \ldots, \tilde{X}_k$. Similarly,
\begin{equation*}
\d \tilde{Y}_2 = |\det (I + \tilde{X}_3 + \cdots + \tilde{X}_k)|^p
  |\det (I + \tilde{X}_2 + \cdots + \tilde{X}_k)|^{-2p} \d \tilde{X}_2
\end{equation*}
and finally,
\begin{equation*}
\d \tilde{Y}_k = |\det (I + \tilde{X}_k)|^{-2p} \d \tilde{X}_k.
\end{equation*}
Since the transformation in (2.3) is of a triangular nature, the
Jacobian matrix will be a triangular block matrix with the Jacobian
being the product of the determinants of the diagonal blocks and the
Jacobian is given by, ignoring the sign,
\begin{align}
\d \tilde{Y}_1 \wedge \ldots \wedge \d \tilde{Y}_k
&=|\det (I + \tilde{X}_1 + \cdots + \tilde{X}_k)|^{-2p}\big\{|\det
(I + \tilde{X}_2 \nonumber\\[.2pc]
&\quad \, + \cdots + \tilde{X}_k)|^{-p} \cdots |\det
(I+\tilde{X}_k)|^{-p}\big\} \d \tilde{X}_1 \wedge \ldots \wedge \d
\tilde{X}_k.
\end{align}
From (2.3), (2.4) and (2.5) we can compute the following product:
\begin{align}
&\bigg\{\prod_{j=1}^k |\det (\tilde{Y}_j)|^{\alpha_j-p}
  |\det (I-\tilde{Y}_j)|^{\alpha_{j + 1} + \cdots + \alpha_{k + 1} + \beta_j + \cdots
  + \beta_k - p}\bigg\} \d \tilde{Y}_1 \wedge \ldots \wedge \d \tilde{Y}_k\nonumber\\
  &\quad\, =\bigg\{\prod_{j = 1}^k|\det (\tilde{X}_j)|^{\alpha_j - p}\bigg \}
  |\det (I + \tilde{X}_2 + \cdots + \tilde{X}_k)|^{\beta_1}\nonumber\\[.2pc]
  &\qquad\ \times\!|\det (I + \tilde{X}_3 + \cdots + \tilde{X}_k)|^{\beta_2}
  \ldots|\det (I + \tilde{X}_k)|^{\beta_{k - 1}}\nonumber\\[.2pc]
  &\qquad\ \times\!|\det (I + \tilde{X}_1 + \cdots + \tilde{X}_k)|^{-(\alpha_1 + \cdots
  +\alpha_{k + 1} + \beta_1 + \cdots + \beta_k)} \d \tilde{X}_1 \wedge \ldots
  \wedge \d \tilde{X}_k.
\end{align}
Multiplying (2.6) on both sides by $c_k$ we have the result since the
right-hand side with $c_k$ is the density in (2.1) and the left-hand side with
$c_k$ is the product of complex matrix-variate type-1 beta densities.

It is easy to see that the converse also holds.
\end{proof}

\begin{theor}[\!]
Let the hermitian positive definite matrices $\tilde{Y}_1, \ldots,
\tilde{Y}_k$ be independently distributed as complex matrix-variate
type-$1$ beta random variables where $\tilde{Y}_j$ has the parameters
$(\alpha_j, \alpha_{j + 1} + \cdots + \alpha_{k + 1} + \beta_j + \cdots + \beta_k)$
for $j = 1, \ldots, k$. Consider the transformation in $(2.3)$ on the space
of $k$-tuples of hermitian positive definite matrices $\tilde{X}_1,
\ldots, \tilde{X}_k$. Then $\tilde{X}_1, \ldots, \tilde{X}_k$ have the
joint density as given in $(2.1)$.
\end{theor}

Thus Theorems~2.1 and 2.2 also provide a characterization of the density
in (2.1). It is known that when $\tilde{Y}_j$ has a complex
matrix-variate type-1 beta density, then $I - \tilde{Y}_j$ again has a
complex matrix-variate type-1 beta density. Thus{\rm ,} from Theorems~2.1 and
2.2 we can get two more results as corollaries. One of them will be
listed here as a theorem and it can also be proved independently by
proceeding parallel to the proof in Theorem~2.1.

\begin{theor}[\!]
Let $\tilde{X}_1, \ldots, \tilde{X}_k$ have the joint density in $(2.1)$.
Consider the transformation
\begin{align}
\tilde{Z}_1 &= (I \!+\! \tilde{X}_1 + \cdots +
\tilde{X}_k)^{-{{1}/{2}}} (I \!+\! \tilde{X}_2 + \cdots +
\tilde{X}_k)(I \!+\! \tilde{X}_1 + \cdots +
\tilde{X}_k)^{-{{1}/{2}}}\nonumber\\[.2pc]
\tilde{Z}_2 &= (I \!+\! \tilde{X}_2 + \cdots +
\tilde{X}_k)^{-{{1}/{2}}}(I \!+\! \tilde{X}_3 + \cdots +
\tilde{X}_k)(I \!+\! \tilde{X}_2 + \cdots + \tilde{X}_k)^{-
{{1}/{2}}}\nonumber\\[.2pc]
& \ \ \vdots\nonumber\\[.2pc]
\tilde{Z}_k &= (I+\tilde{X}_k)^{-1}.
\end{align}
Then $\tilde{Z}_1, \ldots, \tilde{Z}_k$ are independent complex
matrix-variate type-$1$ beta random variables with $\tilde{Z}_j$ having
the parameters $(\alpha_{j + 1} + \cdots + \alpha_{k + 1} + \beta_j + \cdots
+\beta_k,~\alpha_j),$ for $j = 1, \ldots, k$.
\end{theor}

\begin{theor}[\!]
Let $\tilde{X}_1, \ldots, \tilde{X}_k$ have the joint density in
$(2.1)$. Consider the transformation
\begin{align}
\tilde{U}_1 &= \tilde{X}_1^{-{{1}/{2}}}(I + \tilde{X}_2 + \cdots
+ \tilde{X}_k)\tilde{X}_1^{-{{1}/{2}}}\nonumber\\[.2pc]
\tilde{U}_2 &= \tilde{X}_2^{-{{1}/{2}}}(I + \tilde{X}_3 + \cdots + \tilde{X}_k)
 \tilde{X}_2^{-{{1}/{2}}}\nonumber\\[.2pc]
& \ \ \vdots\nonumber\\[.2pc]
\tilde{U}_k &=\tilde{X}_k^{-1}.
\end{align}
Then $\tilde{U}_1, \ldots, \tilde{U}_k$ are independent complex
matrix-variate type-$2$ beta random variables with $\tilde{U}_j$ having
the parameters $(\alpha_{j + 1} + \cdots + \alpha_{k + 1} + \beta_j + \cdots
+ \beta_k,~\alpha_j)$, for $j = 1, \ldots, k$.
\end{theor}

\begin{proof}
From eqs~(2.8), (1.1), (1.2) and Lemma~2.2 we have the following:
\begin{equation*}
\d \tilde{U}_1 = |\det (I + \tilde{X}_2 + \cdots +
\tilde{X}_k)|^p|\det (\tilde{X}_1)|^{-2p} \d \tilde{X}_1
\end{equation*}
for fixed $\tilde{X}_2, \ldots, \tilde{X}_k$.
\begin{equation*}
\d \tilde{U}_2 = |\det (I + \tilde{X}_3 + \cdots + \tilde{X}_k)|^p
|\det (\tilde{X}_2)|^{-2p} \d \tilde{X}_2,
\end{equation*}
and finally
\begin{equation*}
\d \tilde{U}_k = |\det (\tilde{X}_k)|^{-2p} \d \tilde{X}_k.
\end{equation*}
Since the transformation in (2.8) is of a triangular nature, we have the
Jacobian given by
\begin{align}
\!\!\d \tilde{U}_1 \wedge \ldots \wedge \d \tilde{U}_k
&=|\det (\tilde{X}_1)|^{-2p} \ldots |\det (\tilde{X}_k)|^{-2p}
|\det (I \!+\! \tilde{X}_2 + \cdots +\! \tilde{X}_k)|^p\nonumber\\[.2pc]
&\quad \, \times |\det (I + \tilde{X}_3 + \cdots + \tilde{X}_k)|^{p}\ldots|\det
(I + \tilde{X}_k)|^{p}\nonumber\\[.2pc]
&\quad \, \times  \d \tilde{X}_1 \wedge \ldots \wedge \d \tilde{X}_k.
\end{align}
From (2.8),
\begin{align}
I + \tilde{U}_1 &= I + \tilde{X}_1^{-{{1}/{2}}}(I + \tilde{X}_2 + \cdots
 + \tilde{X}_k)\tilde{X}_1^{-{{1}/{2}}}\nonumber\\[.2pc]
 &= \tilde{X}_1^{-{{1}/{2}}}[\tilde{X}_1 + (I + \tilde{X}_2 + \cdots + \tilde{X}_k)]
 \tilde{X}_1^{-{{1}/{2}}}\nonumber\\[.2pc]
 &= \tilde{X}_1^{-{{1}/{2}}}(I + \tilde{X}_1 + \cdots + \tilde{X}_k)\tilde{X}_1^
 {-{{1}/{2}}}\nonumber\\[.2pc]
 I + \tilde{U}_2 &= \tilde{X}_2^{-{{1}/{2}}}(I + \tilde{X}_2 + \cdots + \tilde{X}_k)
 \tilde{X}_2^{-{{1}/{2}}}\nonumber\\[.2pc]
 & \ \ \vdots\nonumber\\[.2pc]
I + \tilde{U}_k &=\tilde{X}_k^{-{{1}/{2}}}(I + \tilde{X}_k)\tilde{X}_k^
{-{{1}/{2}}}.
\end{align}
Now from (2.8), (2.9) and (2.10) we have
\begin{align}
&\bigg\{\prod_{j = 1}^k|\det (\tilde{U}_j)|^{\alpha_{j + 1}
+ \cdots + \alpha_{k + 1} + \beta_j + \cdots + \beta_k - p}\nonumber\\[.2pc]
&\quad\ \times|\det (I + \tilde{U}_j)|^{-(\alpha_j + \cdots + \alpha_{k + 1}
+ \beta_j + \cdots + \beta_k)}\bigg \} \d \tilde{U}_1 \wedge \ldots \wedge
\d \tilde{U}_k\nonumber\\[.2pc]
&=\left \{\prod_{j = 1}^k|\det (\tilde{X}_j)|^{\alpha_j-p}\right \}
|\det (I + \tilde{X}_2 + \cdots + \tilde{X}_k)|^{\beta_1} \cdots |\det (I + \tilde{X}_k)|
^{\beta_{k - 1}}\nonumber\\[.2pc]
&\quad\ \times|\det (I + \tilde{X}_1 + \cdots + \tilde{X}_k)|^{-(\alpha_1 + \cdots +
\alpha_{k + 1} + \beta_k + \cdots + \beta_k)}
\d \tilde{X}_1 \wedge \ldots \wedge \d \tilde{X}_k.
\end{align}
Multiply both sides of (2.11) by $c_k$ to see the result.

The converse also holds. Thus Theorem~2.4 and its converse also provide
a characterization of the density in (2.1). It is known that when
$\tilde{U}_j$ has a complex matrix-variate type-2 beta distribution then
$\tilde{U}_j^{-1}$ has a complex matrix-variate type-2 beta distribution
with the parameters interchanged. This property also gives a couple of
results. Instead of $\tilde{U}_j^{-1}$, a slightly different transformation
will be considered in the next theorem.
\end{proof}

\begin{theor}[\!]
Let $\tilde{X}_1, \ldots, \tilde{X}_k$ have the joint distribution
as given in $(2.1)$. Consider the transformation
\begin{align*}
\tilde{V}_1 &= (I + \tilde{X}_2 + \cdots + \tilde{X}_k)^{-{{1}/{2}}}
 \tilde{X}_1 (I + \tilde{X}_2 + \cdots + \tilde{X}_k)^{-{{1}/{2}}}\\[.2pc]
 \tilde{V}_2 &= (I + \tilde{X}_3 + \cdots + \tilde{X}_k)^{-{{1}/{2}}}\tilde{X}_2
 (I + \tilde{X}_3 + \cdots + \tilde{X}_k)^{-{{1}/{2}}}\\[.2pc]
 & \ \ \vdots\\[.2pc]
 \tilde{V}_{k - 1} &= (I + \tilde{X}_k)^{-{{1}/{2}}}\tilde{X}_{k - 1}(I + \tilde{X}_k)
 ^{-{{1}/{2}}}\\[.2pc]
 \tilde{V}_k &= \tilde{X}_k.
\end{align*}
Then $\tilde{V}_1, \ldots, \tilde{V}_k$ are independent complex
matrix-variate type-$2$ beta random variables with $\tilde{V}_j$ having
the parameters $(\alpha_j,~
\alpha_{j + 1} + \cdots + \alpha_{k + 1} + \beta_j + \cdots + \beta_k)$, for $j = 1,
\ldots, k$.
\end{theor}

The proof can be given by using the steps parallel to the ones in the
proof of Theorem~2.1. The converse of Theorem~2.5 is also true. Further,
Theorem~2.5 and its converse also provide a characterization for the
density in (2.1). By exploiting the relationships between complex
matrix-variate type-1 and type-2 beta random variables one can derive a
number of results and a number of interesting matrix transformations.
These will not be enumerated here in order to save space.

\section*{Acknowledgements}

The last author would like to thank the Natural Sciences and the
Engineering Research Council of Canada for financial assistance. The
authors would like to express their sincere thanks to the referee for
making many valuable suggestions which enabled the authors to make the
presentation far better.


\begin{thebibliography}{99}

\bibitem{1} Biyari~K~H and Lindsey~W~C, Statistical distributions of
hermitian quadratic form in complex Gaussian variables, {\it IEEE Trans.
Information Theory} {\bf 39(3)} (1991) 1076--1082

\bibitem{2} Gupta~R~D and Richards~D~St~P, Multivariate Liouville
distributions, {\it J. Multivariate Anal.} {\bf 23} (1987) 233--256

\bibitem{3} Hayakawa~T, On the distribution of latent roots of a complex
Wishart matrix (non-central case), {\it Ann. Inst.
Stat. Math.} {\bf 24} (1972) 1--17

\bibitem{4} Mathai~A~M, Jacobians of matrix transformations and
functions of matrix argument (New York: World Scientific Publishing)
(1997)

\bibitem{5} Mathai~A~M, An introduction to geometrical probability:
Distributional aspects with applications (New York: Gordon and
Breach Publishers) (1999)

\bibitem{6} Mehta~M~L, Random matrices and statistical theory of energy
levels (New York: Academic Press) (1967)
\end{thebibliography}
\end{document}